\documentclass[12pt]{amsart}
\parskip=\smallskipamount

\newtheorem{theorem}{Theorem}[section]
\newtheorem{lemma}[theorem]{Lemma}
\newtheorem{corollary}[theorem]{Corollary}
\newtheorem{proposition}[theorem]{Proposition}

\theoremstyle{definition}

\newtheorem{example}[theorem]{Example}

\newtheorem{remark}[theorem]{Remark}

\newcommand{\C}{\mathbb{C}}

\newcommand{\R}{\mathbb{R}}

\newcommand{\cC}{\mathcal{C}}

\newcommand{\cO}{\mathcal{O}}

\newcommand\Subset{\subset\!\subset}
\newcommand\wt{\widetilde}

\def\bs{\backslash}

\numberwithin{equation}{section}

\begin{document}
\title[Fibrations and Stein Neighborhoods]
{Fibrations and Stein Neighborhoods}
\author{Franc Forstneri\v c$^*$ \& Erlend Forn\ae ss Wold}
\address{F.\ Forstneri\v c, Institute of Mathematics, Physics and Mechanics, 
University of Ljubljana, Jadranska 19, 1000 Ljubljana, Slovenia}
\email{franc.forstneric@fmf.uni-lj.si}
\address{E.\ F. Wold, Matematisk Institutt, Universitetet i Oslo,
Postboks 1053 Blindern, 0316 Oslo, Norway}
\email{erlendfw@math.uio.no}
\thanks{$^*$Supported by grants P1-0291 and J1-2043-0101, Republic of Slovenia.}

%
%
\date{\today} 

\begin{abstract}
Let $Z$ be a complex space and let $S$ be a compact set in 
$\C^n \times Z$ which is fibered over $\R^n$.
We give a necessary and sufficient condition for $S$ to be a Stein compactum.
\end{abstract}
\maketitle

\section{Fibrations over Totally real sets}
We denote by $\cO(Z)$ the algebra of all holomorphic functions 
on a (reduced, paracompact) 
complex space $Z$, endowed with the compact-open topology. 
A compact subset $K$ of $Z$ is said to be 
a {\em Stein compactum} if $K$ has a basis of Stein open neighborhoods. 
$K$ is $\cO(Z)$-convex if 
\[
		K=\widehat K_{\cO(Z)} = \{z\in Z\colon |f(z)| \le \sup_K |f|,\ \forall f\in\cO(Z)\}.
\]
For the theory of Stein spaces we refer to \cite{GR,Ho}. 

Let $Z$ be a complex space, and consider the product 
$\C^n\times Z$ with the projection $\pi\colon\C^n\times Z\to\C^n$. 
Let $S$ be a compact set in $\C^n \times Z$.  
For every point $\zeta\in\C^n$ we let 
$S_\zeta=\{ z \in Z\colon (\zeta,z)\in S\}$
denote the fiber of $S$ over $\zeta$.
We are interested in the following 

\smallskip
{\em Question:}
Under what conditions on the projection $\pi(S)\subset\C^n$ and on the fibers 
$S_\zeta$ is $S$ a Stein compactum in $\C^n\times Z$?
\smallskip

We give the following precise answer under the
assumption that the projection $\pi(S)$ is 
contained in $\R^n$, the real subspace of $\C^n$.

\begin{theorem}\label{main1}
Let $S$ be a compact set in $\C^n \times Z$ whose projection
$P=\pi(S)$ is contained in $\R^n$.
Then $S$ is a Stein compactum in $\C^n \times Z $ if and only if 
for every open neighborhood $U\subset\mathbb R^n\times Z$ of $S$ 
there exist open sets $V,\Omega\subset \R^n\times Z$,
with $S\subset V\Subset \Omega\subset U$, 
such that for every $u\in P$ the fiber 
$\Omega_u$ is Stein and $\widehat{(S_u)}_{\mathcal O(\Omega_u)}\subset V_u$.  
The analogous result holds if $\pi(S)$ belongs to a 
totally real submanifold of $\C^n$. 
\end{theorem}

A particular reason for looking at this problem is that a question 
of this type, for compact sets that are laminated by holomorphic
leaves, appears in the recent work by the first author \cite{FF:OM};
the relevant result is provided by Corollary \ref{Graphs1} below.

The following simple example (a thin version of Hartog's figure)
shows that it is not enough to assume in Theorem \ref{main1}
that each fiber of $S$ is a Stein compactum.

\begin{example}
Let $Z=\mathbb C$ and let $P \subset \C$ denote the 
real unit interval $P=\{u \in \R \colon 0\leq u \leq 1\}$.
For $0\leq u \leq\frac{1}{2}$ let $S_u =\{z\in\C \colon |z|\leq 1\}$, 
and for $\frac{1}{2}\leq u \leq 1$ let 
$S_u =\{z \in\C \colon \frac{1}{2}\leq |z|\leq 1\}$.  Clearly 
each fiber $S_u$ is a Stein compactum, but due to the continuity 
principle $S$ doesn't have a Stein neighborhood basis in $\C^2$.  
\end{example}

\begin{remark}
In the situation of Theorem \ref{main1}, if $Z$ is a Stein space 
and if every fiber $S_u$ for $u\in \pi(S)\subset\R^n$ 
is $\cO(Z)$-convex, then it is easily seen that $S$ 
is $\cO(\C^n\times Z)$-convex, and hence a Stein compactum.
Following the proof of Proposition 4.3 in \cite{MWO} one 
obtains the following:

\begin{proposition}
Let $\varphi\colon Z\rightarrow X$ be a holomorphic map from a
Stein space $Z$ to a complex manifold $X$.  Let $S\subset Z$
be compact and let $Y=\varphi(S)\subset X$.  If $y\in Y$ is 
a peak point for the algebra $\mathcal O_{X}(Y)$ then 
\[
	\widehat S_{\mathcal O(Z)}\cap S_{y}=\widehat{(S_{y})}_{\mathcal O(Z)}.
\]  
\end{proposition}
\end{remark}

\begin{proof}[Proof of Theorem \ref{main1}]
We shall give the details only when 
$\pi(S) \subset \R^n$; the general case when the projection
is contained in a totally real submanifold 
of $\C^n$ is quite similar.

It is easy to see that the conditions are necessary.  
Indeed, if $S$ is a Stein compactum then for every open
neighborhood $U$ of $S$ in $\R^n\times Z$ there exists a 
Stein neighborhood $\Omega$ of $S$ in $\mathbb C^n\times Z$
such that $\Omega\cap(\R^n \times Z)\subset U$. 
Taking $V\Subset\Omega$ to be an open neighborhood of 
$\widehat S_{\mathcal O(\Omega)}$ in $\R^n\times Z$ we see 
that the hypotheses of the theorem are satisfied. 

Assume now that the conditions hold. 
For $\epsilon>0$ set 
\[
		V(\epsilon) = \{(u+iv,z) \in\C^n\times Z 
		\colon (u,z)\in V,\ \|v\|<\epsilon \}.
\]
To prove the theorem we shall construct
a plurisubharmonic polyhedral neighborhood of $S$ 
contained in $V(\epsilon_0)$ for a given $\epsilon_0>0$. 

Choose an open set $V'$ in $\R^n\times Z$ such that
$V\Subset V'\Subset \Omega$. We shall need two lemmas.

\begin{lemma}
\label{spsh}
There exists a positive strongly plurisubharmonic function
in an open neighborhood $\wt U$ of $S$ in $\C^n \times Z$.
\end{lemma}

\begin{proof}
Fix a point $s_0 =(u_0,z_0)\in S$. 
By compactness it suffices to show that there exists a
plurisubharmonic function in a neighborhood of $S$ 
which is strongly plurisubharmonic near $s_0$.  

Since $V'_{u_0}\Subset \Omega_{u_0}$ and the set $\Omega$ 
is open in $\R^n\times Z$, there exists a $\delta>0$ such that 
$V'_{u} \Subset \Omega_{u_0}$
for all $u$ in the ball 
$B(u_0,2\delta) =\{ \|u-u_0\| < 2\delta\}$.  

As $W = (B(u,2\delta) + \mathrm{i} \R^n) \times \Omega_{u_0}$ 
is Stein, there exists a positive strongly plurisubharmonic 
function $\rho_{1}$ on $W$.  We may assume that $\rho_1(s_0)=1$.  

Choose $M\in\R$ such that for any $u \in\partial B(u_0,\delta)$
and any $z\in V'_{u}$ we have that $\rho_{1}(u,z)<M$.  
Let $\rho_{2}$ be a plurisubharmonic function on $\C^n$ 
such that $\rho_{2}(u_0)=0$ and $\rho_{2}(u)>M$
for all $u\in\partial B(u_0,\delta)$. 
If $\epsilon>0$ is chosen small enough then the function  
$\rho=\mathrm{max}\{\rho_{1},\rho_{2}\}$ is well defined 
and plurisubharmonic on $V'(\epsilon)$, and is 
strongly plurisubharmonic near the point $s_0$. 
\end{proof}

By shrinking the set $U\supset S$ we may assume that 
$U=\wt U\cap (\R^n \times Z)$, where $\wt U$ satisfies
the conclusion of Lemma \ref{spsh}.

\begin{lemma}
\label{sepfunction}
For any point $q$ belonging to the boundary of $V$ 
in $\R^n\times Z$ there exists an $\epsilon>0$ and
a continuous plurisubharmonic function $\rho$ on 
$V'({\epsilon})$ such that $\rho(q)=2$ and $\rho <\frac{1}{2}$ on $S$.
\end{lemma}

\begin{proof}
If $\pi(q)\notin\pi(S)=P$, we may find a holomorphic function 
$g\in\mathcal O(\C^n)$ such that $\rho=|g\circ\pi|$ will work.  

Assume now that $q=(u,z)$ with $u \in P$. 
Since $\widehat{(S_u)}_{\mathcal O(\Omega_u)}\subset V_u$
by the assumption, there is a holomorphic function $f\in\cO(\Omega_u)$ 
such that $f(z)=2$ and $|f|<\frac{1}{2}$ on $S_u$.  
By continuity there is a $\delta>0$ 
such that for all $u'\in \overline{B(u,\delta)}$ we have that
$V'_{u'}\Subset \Omega_u$ and $|f|<\frac{1}{2}$ on $S_{u'}$ .   
Let $M\in\R$ be such that $|f| \leq M$ on $V'_{u'}$ 
for all $u'\in \overline{B(u,\delta)}$.
We consider $f$ as a holomorphic function on 
$\C^n\times \Omega_u$ which is independent of the first variable.

Let $\chi\in\mathcal C^\infty_0(B(u,\delta))$ be a smooth function with compact
support such that $0\leq\chi\leq 1$ and $\chi(u)=1$.  
Let $g\in\cO(\C^n)$ be a holomorphic function approximating 
$\chi$ close enough on $\overline{B(u,\delta)}$
such that $\|g\|_{\partial B(u,\delta)}<\frac{1}{3M}$
and $g(u)=1$. Then $\rho_1:=|f\cdot g|$ is plurisubharmonic on 
$[\overline{B(u,\delta)}\oplus\mathrm{i\mathbb R^n}]\times \Omega_u$, 
$\rho_1(q)=|f(z)|=2$, 
and $\rho_1(w)<\frac{1}{3}$ for all $w=(u',z)$ such that 
$u'\in\partial B(u,\delta)$ and $z \in V'_{u'}$. 

Let $h\in\mathcal O(\mathbb C^n)$ be a holomorphic function such that 
$h(u)=0$, $0\leq |h|\leq\frac{1}{2}$ on $B(u,\delta)$, 
$\frac{1}{3}<|h|<\frac{1}{2}$ on $\mathbb R^n\bs B(u,\delta)$.   
Such $h$ exists since continuous functions on $\mathbb R^n$ 
can be approximated uniformly by entire functions on $\C^n$
(see for instance \cite{Sc}).  

The function $\rho_2:=|h\circ \pi|$ is then 
plurisubharmonic on $\C^n \times Z$.   

If $\epsilon>0$ is chosen small enough then the function 
$\rho=\mathrm{max}\{\rho_1,\rho_2\}$ is well defined on $V'(\epsilon)$
and it satisfies the conclusion of Lemma \ref{sepfunction}.
\end{proof}

We can now complete the proof of Theorem \ref{main1}.
By compactness and Lemma \ref{sepfunction} there exist an 
$\epsilon>0$, with $2\epsilon <\epsilon_0$, 
and finitely many plurisubharmonic functions $\rho_1,\ldots,\rho_m$ 
on $V'(2\epsilon)$ such that $\rho_j<\frac{1}{2}$ on $S$ for $j=1,\ldots,m$,
and such that for every $q=(u,z)\in \partial V \subset \R^n\times Z$
we have $\rho_j(q)>\frac{3}{2}$ for at least one $j\in\{1,\ldots,m\}$.

Denote the variable on $\C^n$ by $\zeta=u+iv$.
If we replace $\rho_j$ by $\rho_j+C|v|^2$ for a sufficiently large 
$C>0$, then we have that 
\[
	\partial [V(\epsilon)] \subset 
	\bigcup_{j=1}^m \, \left\{w\in V'(2\epsilon) \colon 
	\rho_j(w) > \frac{3}{2}\right\}.
\]  
Let $W$ denote the set 
\[
		W = \{w\in V'(2\epsilon)\colon \rho_j(w)<1,\ j=1,...,m\}.  
\]    
Let $W_0$ be the union of all connected components of $W$ 
which intersect $S$. Then $S\subset W_{0}$ and 
$W_0$ is relatively compact in $V(\epsilon)$.   

Let $\varphi\in\cC^\infty((-\infty,1))$ be a convex increasing function 
such that $\varphi(t)\rightarrow +\infty$ as $t\rightarrow 1$.   Let 
$\wt \rho_j=\varphi\circ\rho_j$ for $j=1,...,m$.  
Then 
\[
		\rho:=  
             \max_{1\le j\le m} \wt\rho_j
\]
is a plurisubharmonic exhaustion function of $W_0$.  
By Lemma \ref{spsh} we may add to $\rho$ a strongly plurisubharmonic function 
and get a strongly plurisubharmonic exhaustion function
of $W_0$, and hence $W_0$ is Stein according to Narasimhan's theorem \cite{Na}.
Since $\overline W_0\subset V(\epsilon_0)$, this concludes the proof. 
\end{proof}

\begin{corollary}
\label{Cor1}
Assume that $S$ is a compact set in $\R^n\times Z$ 
such that for any open neighborhood $U\subset \R^n\times Z$ 
of $S$ and for each $u\in \pi(S)\subset\R^n$ 
there exist a Stein neighborhood $\Omega_u\Subset U_u$
of the fiber $S_u$ and a number $\delta>0$ such that $S_{u'}$ is 
holomorphically convex in $\Omega_u$ for each $u'$ 
with $\|u-u'\|\leq\delta$. 
Then $S$ is a Stein compactum in $\C^n\times Z$. 
\end{corollary}

\begin{proof}
By compactness there is a finite number of $u_j$'s, 
$\Omega_{u_j}$'s and $\delta_j$'s such that 
$\{B(u_j,\delta_j)\}$ is an open cover of $\pi(S)$ in $\R^n$. 
If all $\delta_j$'s are small enough 
we get that $\bigcup_j B(u_j,\delta_j)\times\Omega_{u_j}$ is contained in $U$. 

Define a neighborhood $\Omega$ of $S$ in $\R^n\times Z$ as follows.  
For a point $u\in\overline{B(u_j,\delta_j)}$ let 
the fiber $\Omega_u$ consist of the intersection of all $\Omega_{u_k}$ such that 
$u$ is also contained in $\overline{B(u_k,\delta_k)}$.   Then $\Omega$ is an open 
neighborhood of $S$ contained in $U$, and $S_u$ is holomorphically convex 
in $\Omega_u$ for every $u$.  Let $V\Subset\Omega$ be any open neighborhood 
that contains $S$.  The conditions in Theorem \ref{main1} are satisfied
and hence the corollary follows. 
\end{proof}

\section{Stein neighborhoods of certain laminated sets} 
According to Siu's theorem \cite{Si}
every closed Stein subvariety of a complex space admits
an open Stein neighborhood in that space. 
For simpler proofs and generalizations
to $q$-convex subspaces see \cite{Co,De}.

The following generalization of Siu's theorem is 
proved in \cite[Theorem 2.1]{FF:EOP}
(for the last stetement see also \cite[Theorem 1.2]{FF:Kohn}).

\begin{theorem}
\label{FF-EOP}
Let $X$ be a closed Stein subvariety of complex space $Z$.
Assume that $S$ is a compact set in $Z$ that is
$\cO(\Omega)$-convex in an open Stein domain $\Omega\subset Z$
containing $S$ and such that $S \cap X$ is $\cO(X)$-convex.
Then for every open set $U$ in $Z$ containing $S\cup X$
there exists an open Stein domain $V$ in $Z$ such that 
$S\cup X \subset  V\subset U$ and $S$ is $\cO(V)$-convex.
\end{theorem}

In the remainder of the article we consider the following situation: 

$Z$ is a complex space,
$X$ is a Stein space, $h\colon Z\to X$ is a holomorphic submersion 
onto $X$, $K$ is a compact $\cO(X)$-convex subset of $X$,  and
$P$ is a compact set in $\R^n$ (considered as the real subspace
of $\C^n$). 

The following corollary is used in an essential way
in the proof of the main result in \cite{FF:OM}
to the effect that all Oka properties of a complex manifold
are equivalent to each other.

\begin{corollary}
\label{Graphs1}
{\rm (Assumptions as above.)} 
Let $f\colon P\times X\to Z$ be a continuous map such that 
$h\circ f(p,x)=x$ for all $(p,x)\in P\times X$, 
and such that $f_p=f(p,\cdotp)\colon X \to Z$ is holomorphic for 
every $p\in P$. Then the set  
\begin{equation}
\label{Sigma}
	\Sigma =\{(p,f(p,x)) \colon p\in P,\ x\in K\} 
\end{equation}
admits an open Stein neighborhood $\Theta$ in $\C^n \times Z$
such that $\Sigma$ is $\cO(\Theta)$-convex.
\end{corollary}

\begin{proof}
For every $p\in P$ the set $V_p=f_p(X)$ is a closed Stein
subvariety of $Z$, and $\Sigma_p = f_p(K) \subset V_p$ 
(the fiber of $\Sigma$ over $p$) is $\cO(V_p)$-convex.

Fix $p\in P$. By Siu's theorem \cite{Si} there exists an open 
Stein neighborhood $\Omega \subset Z$ of $\Sigma_p$.
If $q\in P$ is sufficiently near $p$ then, due to the continuity of $f$,
we have that $\Sigma_q \subset V_q\cap \Omega$.
Since $\Sigma_q$ is $\cO(V_q)$-convex, it is also 
$\cO(V_q\cap \Omega)$-convex.
As $V_q\cap \Omega$ is a closed subvariety of the Stein domain
$\Omega$, it follows that $\Sigma_q$ is also $\cO(\Omega)$-convex.
(Indeed, for any point $z\in \Omega\bs V_q$ there exists by
Cartan's theorem a holomorphic function on $\Omega$ 
that equals one at $z$ and that vanishes on $V_q\cap \Omega$;
hence no such point can belong to the $\cO(\Omega)$-hull of $\Sigma_q$.)

This shows that $\Sigma$
satisfies the assumptions of Corollary \ref{Cor1} and hence 
it admits a basis of open Stein neighborhoods in $\C^n \times Z$.

Since $P$ is contained in $\R^n$, we have that 
$\mathcal O(\mathbb C^n)|_P$ is dense in $\mathcal C(P)$,
the space of complex valued continuous functions on $P$ (see \cite{Sc}). 
By using cutoff functions in the Euclidean variable
and approximating them sufficiently well by entire 
functions  we see  that $\Sigma$ is $\cO(\Theta)$-convex 
in every open Stein set $\Theta\subset \C^n \times Z$
containing $\Sigma$. This completes the proof.
\end{proof}

\begin{remark}
The set $\Sigma$ in Corollary \ref{Graphs1}
is laminated by the graphs of the holomorphic functions 
$f_p$ over $K$, depending continuously on
the parameter $p\in P$. This lamination is rather simple since 
it admits a global product structure.
In general, the existence of Stein neighborhoods of sets 
that are foliated (or laminated) by complex submanifolds
is a rather subtle issue as is seen in the case of
worm type Levi-flat hypersurfaces (see \cite{BeF,DF,FCL}).
Corollary \ref{Graphs1} extends
to the more general situation when the leaves $\Sigma_p$ 
are not necessarily graphs, but compact holomorphically convex 
subsets in a continuously moving family of 
Stein subvarieties of $Z$.
\qed \end{remark}

In the final result we shall need the following
well known implication of Rossi's local maximum modulus principle
(see e.g.\ \cite[Lemma 6.5]{FF:INT} for the 
case $\Omega=\C^N$; the general case easily reduces to this one 
since we can replace $\Omega$ by a  relatively compact
Stein subset which then embeds as a closed complex subvariety
in a Euclidean space). 

\begin{lemma}
\label{Rossi}
Assume that $V$ is a closed complex subvariety of a Stein space $\Omega$.
If $S$ is a compact $\cO(\Omega)$-convex subset of $\Omega$
and if $K$ is a compact $\cO(V)$-convex subset of $V$ 
such that $S\cap V \subset K$, then $K\cup S$ is $\cO(\Omega)$-convex.
\end{lemma}

Our final result furnishes a parametric version of Theorem \ref{FF-EOP}.

\begin{corollary}
\label{Graphs2}
Assume the situation of Corollary \ref{Graphs1}.
Set $V_p=f_p(X)\subset Z$ for $p\in P$, and let $\Sigma \subset P\times Z$
be the compact set  (\ref{Sigma}) with the fibers $\Sigma_p=f_p(K) \subset V_p$.
Let $S$ be a compact set in $P\times Z \subset \C^n\times Z$, 
with fibers $S_p$ $(p\in P)$, such that
\begin{itemize}
\item[(i)]  $S_p\cap V_p \subset \Sigma_p$ and $S_p\cap V_p$ 
is $\cO(V_p)$-convex for each $p\in P$, and
\item[(ii)]  $S$ is $\cO(\Theta)$-convex in an open Stein domain  
$\Theta\subset \C^n \times Z$.  
\end{itemize}
Then $\Sigma \cup S$ is a Stein compact in $\C^n \times Z$.
\end{corollary}

\begin{proof}
Fix $p\in P$. By Theorem \ref{FF-EOP} the set $(\{p\}\times V_p) \cup S$ 
admits an open Stein neighborhood $\Omega\subset \C^n \times Z$ such 
that $S$ is $\cO(\Omega)$-convex. 
It follows that the fiber $S_q$ is $\cO(\Omega_q)$-convex for every $q\in P$.
(Here $\Omega_q$ denotes the fiber of $\Omega$ over the point $q\in P$.)

Choose an open Stein neighborhood $\Theta\subset Z$ of the 
compact set $\Sigma_p\cup S_p$ such that 
$\overline \Theta$ is compact and contained in $\Omega_p$.
Since $S$ is compact, $\Omega$ is open and $f$ is continuous,
we have for all $q\in P$ sufficiently near $p$ that
\[
		\Sigma_q\cup S_q\subset \Theta \subset \Omega_q.
\]
For such $q$ we see as in the proof of Corollary \ref{Graphs1}
that $\Sigma_q$ is $\cO(\Theta)$-convex. 
Since $S_q$ is also $\cO(\Omega_q)$-convex and hence 
$\cO(\Theta)$-convex, Lemma \ref{Rossi} now shows that 
$\Sigma_q \cup S_q$ is $\cO(\Theta)$-convex for all
$q\in P$ sufficiently near $p$. The conclusion now follows from
Corollary \ref{Cor1}.
\end{proof}

\bibliographystyle{amsplain}

\begin{thebibliography}{10}



\bibitem{BeF}
BEDFORD, E., FORN\AE SS, J.\ E.,
Domains with pseudoconvex neighborhood systems.
Invent.\ Math., 47 (1978), 1--27.

\bibitem{Co}
COL\c TOIU, M., 
Complete locally pluripolar sets. 
J.\ Reine Angew.\ Math., 412 (1990), 108--112. 


\bibitem{De} 
DEMAILLY, J.-P.,
Cohomology of $q$-convex spaces in top degrees.
Math.\ Z., 204 (1990), 283--295.

\bibitem{DF} 
DIEDERICH, K., FORN\AE SS, J.-E.,
Pseudoconvex domains: an example with nontrivial Nebenh\"ulle. 
Math.\ Ann., 225 (1977), 275--292.

\bibitem{FF:INT} 
FORSTNERI\v C, F.,
Interpolation by holomorphic automorphisms and
embeddings in $\C^n$.
J.\ Geom.\ Anal., 9 (1999), 93--118.

\bibitem{FF:EOP}
FORSTNERI\v C, F.,
Extending holomorphic mappings from subvarieties in Stein manifolds.
Ann.\ Inst.\ Fourier, 55 (2005), 733--751.


\bibitem{FF:Kohn}
FORSTNERI\v C, F.,
The Oka principle for sections of stratified fiber bundles.
Pure Appl.\ Math.\ Quarterly 
(Special Issue in honor of Joseph J. Kohn), 
6 (2010), no.\ 3, 843--874.

\bibitem{FF:OM}
FORSTNERI\v C, F.,
Oka manifolds.
C.\ R.\ Acad.\ Sci.\ Paris., Ser.\ I 347 (2009), 1017--1020.

\bibitem{FCL} 
FORSTNERI\v C, F., LAURENT-THI\'EBAUT, C.,
Stein compacts in Levi-flat hypersurfaces 
Trans.\ Amer.\ Math.\ Soc., 360 (2008), no.\ 1, 307--329. 

\bibitem{GR} 
GUNNING, R.\ C., ROSSI, H., 
\textit{Analytic functions of several complex variables.}
Prentice--Hall, Englewood Cliffs, 1965.

\bibitem{Ho} 
H\"ORMANDER, L., 
\textit{An Introduction to Complex Analysis in Several Variables}. 
Third ed. 
North Holland, Amsterdam, 1990.

\bibitem{MWO}
MANNE, P.\ E., WOLD, E.\ F., \O VRELID, N.,
Carleman approximation by entire functions on 
Stein manifolds.  Preprint, University of Oslo (2008).  

\bibitem{Na}
NARASIMHAN, R.,
The Levi problem for complex spaces.
Math.\ Ann., 142 (1961), 355--365.

\bibitem{Sc}
SCHEINBERG, S.,
Uniform approximation by entire functions.
J.\ Analyse.\ Math., 29 (1976), 16--18.

\bibitem{Si}
SIU, J.-T.,
Every Stein subvariety admits a Stein neighborhood.
Invent.\ Math., 38 (1976), 89--100.


\end{thebibliography}

\end{document}